\documentclass{amsart}
\usepackage{graphicx}
\usepackage{fullpage}
\usepackage{cancel}
\usepackage{epstopdf}
\usepackage{tikz}
\usepackage{amsmath}
\usepackage{microtype}
\usepackage[numbers]{natbib}
\usepackage{algorithm,algorithmic}
\usepackage{amsaddr}

\author{Sumedh Joshi}
\address[Sumedh Joshi]{Center for Applied Mathematics \\ 657 Rhodes Hall, Cornell University \\ Ithaca NY, 14850}

\author{Peter Diamessis}
\address[Peter Diamessis]{School of Civil and Environmental Engineering \\ 105 Hollister Hall, Cornell University \\ Ithaca NY, 14850}

\usepackage{amsmath}   
\usepackage{amsthm}
\usepackage{amssymb} 
\usepackage{calc}
\usepackage{ifthen}
\usepackage{graphicx}
\usepackage{enumerate}
\usepackage{centernot}
\usepackage[mathscr]{eucal}
\renewcommand\and{\ensuremath{\;\;\;\;\;\;\text{and}\;\;\;\;\;\;} }

\newcommand\restr[2]{{% we make the whole thing an ordinary symbol
  \left.\kern-\nulldelimiterspace % automatically resize the bar with \right
  #1 % the function
  \vphantom{\big|} % pretend it's a little taller at normal size
  \right|_{#2} % this is the delimiter
}}

 \newcommand\R{\ensuremath{\bb{R}}}

 \newcommand\N{\ensuremath{\bb{N}}}

\newcommand\norm[1]{\ensuremath{\left\lvert \left\lvert  #1 \right\rvert \right\rvert }}

\newcommand\avg[2][1]{
	\ifthenelse{\equal{#1}{1}}{\ensuremath{\left \langle#2\right \rangle }}{}
	\ifthenelse{\equal{#1}{2}}{\ensuremath{\left \langle#2^2\right \rangle }}{}
	\ifthenelse{\equal{#1}{3}}{\ensuremath{\left \langle#2\right \rangle ^2}}{}}

\newcommand\derv[3][1]{
	\ifthenelse{\equal{#1}{1}}{\ensuremath{\frac{d#2}{d#3}}}{}
	\ifthenelse{\equal{#1}{2}}{\ensuremath{\frac{d^2#2}{d#3^2}}}{}}

\newcommand\prtl[3][1]{	\ifthenelse{\equal{#1}{1}}{\ensuremath{\frac{\partial#2}{\partial#3}}}{}
	\ifthenelse{\equal{#1}{2}}{\ensuremath{\frac{\partial^2#2}{\partial#3^2}}}{}}

\newcommand\bb[1]{\ensuremath{\mathbb{#1}}}

\renewcommand\eqref[1]{Eq. (\ref{#1})}

\theoremstyle{plain}   
\theoremstyle{definition}  
\numberwithin{exer}{subsection} % important bit
\theoremstyle{plain}   
\theoremstyle{plain}   
\theoremstyle{definition}   
\theoremstyle{definition}   
\theoremstyle{plain}   
\theoremstyle{plain}   
\theoremstyle{plain}    
\title{A deflated Schur complement method for the iterative solution of a high-order discontinuous element discretization of the Poisson equation}

\begin{document}
\maketitle

\begin{abstract}
A combination of block-Jacobi and deflation preconditioning is used to solve a high-order discontinuous element-based collocation discretization of the Schur complement of the Poisson-Neumann system as arises in the operator splitting of the incompressible Navier-Stokes equations.
The ill-posedness of the Poisson-Neumann system manifests as an inconsistency of the 
Schur complement problem, but it is shown that this can be accounted for with appropriate projections out of the null space of the Schur complement matrix without affecting the accuracy of the solution.  The block-Jacobi preconditioner, combined with deflation, is shown to yield GMRES convergence independent of the polynomial order of expansion within an element.  Finally, while the number of GMRES iterations does grow as the element size is reduced (e.g. $h$-refinement), the dependence is very mild; the number of GMRES iterations roughly doubles as the element size is divided by a factor of six.  In light of these numerical results, the deflated Schur complement approach seems practicable, especially for high-order methods given its convergence independent of polynomial order.
\end{abstract}

\section{Introduction}

\subsection{Background}

Domain decomposition methods have been recently applied to high-order, discontinuous, discretizations of elliptic problems with good success \cite{Antonietti2010,Antonietti2007,Brix2015,Canuto2013,Hesthaven}.  These methods aim to build preconditioning strategies for iterative Krylov-type solvers to obtain an accurate solution whose convergence is independent of the parameters of the grid.  These approaches have largely focused on the widely-used discontinuous Galerkin class of numerical methods, and many leverage the symmetric positive-definite nature of the DG discretization.

% Paragraph: What is the Schur complement problem?

Analogous to the DG class of discretization methods, the Spectral Multidomain Penalty Method (SMPM) is a high-order discontinuous variant of the spectral element method that uses spectral differentiation matrices to compute derivatives \cite{Hesthaven1998}.  Because spectral differentiation matrices are themselves unsymmetric \cite{Costa2000} the operator matrices resulting from SMPM are unsymmetric and not self-adjoint. Nevertheless, the SMPM discretization has been used to solve complex, large-scale, environmental fluid mechanics problems on hundreds of processors and with hundreds of millions of unknowns \cite{Diamessis2006,Abdilghanie2013,Diamessis2011}.

\subsection{Spectral multi-domain penalty method}

In this work, we describe a method for solving the 2D Poisson-Neumann system that arises within the time-splitting of the 2D incompressible Navier-Stokes equations \cite{Escobar-Vargas2014}, and is given on a domain $\Omega \subset \R^2$ as
	\begin{align}
		\nabla^2 u &= f  \textrm{    on    } \Omega \nonumber \\
		n\cdot \nabla u &= g \textrm{    on   } \partial\Omega.
		\label{poisson}
	\end{align}
$\Omega$ is discretized into an $m_y \times m_x$ cartesian quadrilateral element grid with elements $\Omega_{ij}$, where $i = \{1, \dots, m_y\}$ and $j = \{1,\dots, m_x\}$.  Within each element is a 2D Guass-Lobatto-Legendre (GLL) grid with $n$ GLL points per direction for a total of $n^2$ grid points per element.  Since the discretization is discontinuous, function values are allowed to differ along the $2n$ grid points on the boundary of each pair of elements, and thus the full grid has a total of $r = n^2 m_x m_y$ grid points.

Now we define the SMPM element matrices and inter-element continuity conditions. 
   Let $Lu = f$ represent the discrete Poisson-Neumann system on $\Omega \subset \R^2$ a domain discretized into an $m_x \times m_y$ element mesh with each element $V_i$ smoothly and invertibly mapped from the master element $[-1,1]\times[-1,1]$.
   On each element a two-dimensional Gauss-Lobatto-Legendre (GLL) grid with $n$ points in each direction is constructed and used to evaluate the Lagrange interpolant basis and their derivatives by way of spectral differentiation matrices\cite{Costa2000}.  Thus each element contains $n^2$ grid points.
   If $V_i$ and $V_j$ share the $n$ GLL points along one of their  four boundaries, then each element owns a copy of those $n$ GLL nodes in order maintain the discontinuous nature of this method.
   Thus as a matrix, $L \in \R^{r\times r}$ is of dimension $r = n^2 m_xm_y$, where $r$ denotes the total number of nodes in the grid $\Omega$.
      
   In the SMPM the weak inter-element continuity condition is of Robin type, and is enforced by the flux $R_{ij} : \partial V_j \longrightarrow \partial V_i$ from element $V_j$ into $V_i$ for $V_i,V_j$ with an adjacent boundary $\partial V_j \cap \partial V_i$ consisting of $n$ grid points.  $R_{ij}$ is defined as
   \begin{align}
      R_{ij} = I + \hat{n}_i \cdot \nabla
   \end{align}
   where $\hat{n}_i: \partial V_i \longrightarrow \R^2$ is the outward pointing normal vector of $\partial V_i$ and $I$ is the identity operator.  A depiction of a $2 \times 2$ element grid with the inter-element fluxes is shown in Fig. 1, in which the elements $V_1,V_2,V_3,V_4$ have been separated to emphasize the discontinuous nature of the SMPM.
   
\begin{figure}
	\begin{tikzpicture}
   \label{fig:2x2}

		% The elements.
		\draw[line width=0.25mm](0,0) rectangle(2,2);
		\draw[line width=0.25mm](3,3) rectangle(5,5);
		\draw[line width=0.25mm](0,3) rectangle(2,5);
		\draw[line width=0.25mm](3,0) rectangle(5,2);

		% Draw the GLL grid in each of the four elements.
		\foreach \x in { -1.0, -0.93, -0.78, -0.57, -0.30, 0, 0.30, 0.57, 0.78, 0.93, 1.0 }
			\foreach \y in { -1.0, -0.93, -0.78, -0.57, -0.30, 0, 0.30, 0.57, 0.78, 0.93, 1.0 }
				\draw node[fill,circle,scale=0.1](\x,\y) at (\x+1,\y+1) {};

		\foreach \x in { -1.0, -0.93, -0.78, -0.57, -0.30, 0, 0.30, 0.57, 0.78, 0.93, 1.0 }
			\foreach \y in { -1.0, -0.93, -0.78, -0.57, -0.30, 0, 0.30, 0.57, 0.78, 0.93, 1.0 }
				\draw node[fill,circle,scale=0.1](\x,\y) at (\x+4,\y+1) {};

		\foreach \x in { -1.0, -0.93, -0.78, -0.57, -0.30, 0, 0.30, 0.57, 0.78, 0.93, 1.0 }
			\foreach \y in { -1.0, -0.93, -0.78, -0.57, -0.30, 0, 0.30, 0.57, 0.78, 0.93, 1.0 }
				\draw node[fill,circle,scale=0.1](\x,\y) at (\x+4,\y+4) {};

		\foreach \x in { -1.0, -0.93, -0.78, -0.57, -0.30, 0, 0.30, 0.57, 0.78, 0.93, 1.0 }
			\foreach \y in { -1.0, -0.93, -0.78, -0.57, -0.30, 0, 0.30, 0.57, 0.78, 0.93, 1.0 }
				\draw node[fill,circle,scale=0.1](\x,\y) at (\x+1,\y+4) {};

		% Denote the subdomains.
		%\node at (1,5.5) {{$\Omega_1 = V_1 \cup V_2$}};
		%\node at (4,5.5) {{$\Omega_2 = V_3 \cup V_4$}};

		% The labels for the elements.
		\node at (1,1) {\huge{$V_2$}};
		\node at (4,4) {\huge{$V_3$}};
		\node at (1,4) {\huge{$V_1$}};
		\node at (4,1) {\huge{$V_4$}};

		% The arrows between the elements to show inter-element fluxes.
		\draw [line width=0.5mm, ->] (2,1.5) -- node[above]{\small{$R_{42}$}}  (3,1.5);
		\draw [line width=0.5mm, ->] (3,0.5) -- node[below]{\small{$R_{24}$}} (2,0.5);

		\draw [line width=0.5mm, ->] (2,4.5) -- node[above]{\small{$R_{31}$}}  (3,4.5);
		\draw [line width=0.5mm, ->] (3,3.5) -- node[below]{\small{$R_{13}$}} (2,3.5);

		\draw [line width=0.5mm, ->] (0.5,2) -- node[left]{\small{$R_{12}$}}  (0.5,3);
		\draw [line width=0.5mm, ->] (1.5,3) -- node[right]{\small{$R_{21}$}} (1.5,2);

		\draw [line width=0.5mm, ->] (3.5,2) -- node[left]{\small{$R_{34}$}}  (3.5,3);
		\draw [line width=0.5mm, ->] (4.5,3) -- node[right]{\small{$R_{43}$}} (4.5,2);

	\end{tikzpicture}
	\caption{A depiction of the logical arrangement of a  $2\times 2$ element spectral multi-domain penalty method (SMPM) grid with $10 \times 10$ Gauss-Lobatto-Legendre points in each element denoted $V_j$ for $j = 1,2,3,4$. The inter-element continuity fluxes are represented with $R_{ij}$ with $i\neq j$. }
\end{figure}
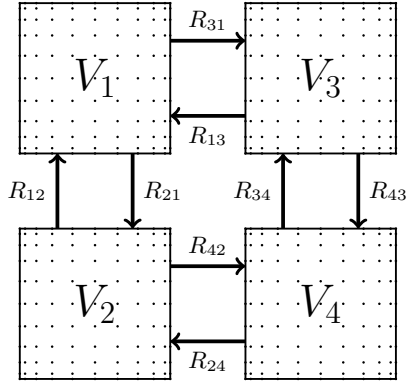

The physical boundary conditions are Neumann, and are given on $\partial V_i \cap \partial \Omega$ as $n_i \cdot \nabla$ where $n_i$ is again the outward pointing normal vector.
   Given a function $u$, on an element $V_i$ the residual in the spectral multi domain penalty method is given by the sum of the Laplacian, the inter-element continuity mismatch, and the boundary condition mismatch as
	\begin{align}
      \label{smpm_poisson}
      L_i u_i = \nabla^2 u_i + \tau_i\left( R_{ii} u_i  - \sum_{j\in N(i)} \restr{R_{ij} u_j}{\partial V_i \cap \partial V_j}  \right) + \tau_i \restr{\hat{n}_i \cdot \nabla u_i}{\partial V_i \cap \partial \Omega}  = f_i + \tau_i g_i.
	\end{align}
   Here, $g_i$ is the boundary value of the Neumann boundary condition restricted to element $V_i$, and $N(i)$ is the index set of elements adjacent to $V_i$.  The inter-element continuity, external boundary conditions, and the PDE are all satisfied weakly, since the residual is the sum of these three components.  The penalty parameter $\tau_i$ represents the degree to which the inter-element continuity and boundary conditions are weighted in the residual relative to the PDE, and the optimal choice of $\tau_i$ is determined by stability criteria for hyperbolic problems\cite{Hesthaven1997,Hesthaven1998}, and a heuristic for the Poisson problem \cite{Escobar-Vargas2014}.

\subsection{Construction of the Schur complement problem}

% The point: What the domain decomposition we perform looks like, coarsely speaking.
As shown in Figure \ref{fig_schur_preconditioner_grid} the domain $\Omega$ is discretized with a collection of \emph{elements} $V_j$, each invertibly mapped from the master element $[-1,1]\times [-1,1]$.  $\Omega$ is decomposed into $m_xm_y$ many \emph{sub-domains} $\Omega_i$, with each sub-domain corresponding to a single element in the mesh.  Since the SMPM is a high-order method, each element represents large, local, and dense linear algebraic operations.  Thus, it makes sense from a domain decomposition perspective for the elements and sub-domains to coincide.

 Along each of the $(m_x - 1)(m_y - 1)$ interfaces between the sub-domains/elements are $2n$ GLL nodes ($n$ nodes on either side of each interface).
   Denote as $k$ the number of interfacial nodes in the domain decomposition, and this set of $k$ interface nodes as $\Gamma$.
   The discrete Poisson operator $L$ (Eq.~\ref{smpm_poisson}) is decomposed into a local term and an inter-subdomain flux term which is used to construct the Schur problem.
   This operator decomposition comprises three operators which are defined below.

   First, denote as $E:\Gamma \longrightarrow \Omega$ the inclusion map that maps from the interfacial grid $\Gamma$ to $\Omega$.
   As a matrix, $E\in \R^{r \times k}$ and is composed of zeros and ones, $E^T$ is the restriction from the full grid to the interface grid $\Gamma$, and $E^T E = I \in \R^{k \times k}$ the identity matrix. Naturally $EE^T$ is not an identity matrix.

Second, define an operator $B : \Omega \longrightarrow \Gamma$ that consists of the inter-subdomain Robin boundary fluxes.
$B$ represents all of the the inter-element fluxes $R_{ij}$.
As a matrix, $B \in \R^{k \times r}$, since it computes $I + \hat{n} \cdot \nabla$ within a subdomain using spectral differentiation matrices and assigns it to the interface of its neighbor.

Finally, define the operator $A : \Omega \longrightarrow \Omega$, which represents the part of $L$ that is entirely local to one subdomain.  $A$ consists of the Laplacian part of $L$, the boundary condition mismatch,  and the $R_{ii}$ terms in Eq.~(\ref{smpm_poisson}).
Since $A$ is entirely local to each subdomain, as a matrix $A\in \R^{r\times r}$ is block-diagonal.  $A$ represents $m_xm_y$ decoupled homogenous Poisson-Robin boundary value problems, and as such is invertible and block diagonal.

These three operators are defined so that their combination yields the SMPM Poisson-Neumann operator
\begin{align} 
   Lu = Au + EBu = f.
\end{align}
Notice that the action of $EB$ couples the subdomains only weakly since its action is combined with $A$ in the residual $(A + EB)u = f$.  This weak enforcement of inter-subdomain continuity in the SMPM allows for decoupling the subdomains by decoupling the action of $B$ from that of $A$.  To accomplish this, we seek a vector $v \in \R^k $ on the interfacial nodes, $\bigcap_i \partial \Omega_i$, such that the solution to $A_i u_i = f_i - Ev_i$ on each $\Omega_i$ subdomain also solves $Lu = f$.  By writing 
\begin{align}
	Au = f - EBu
\end{align}
it is clear that $v = Bu$, the image of the solution under the inter-subdomain flux operator.  Because $B$ is a contraction (mapping from the full grid to the interfacial grid), finding the image $Bu$ is easier than finding $u$ itself; its value is given by the solution to the system
\begin{align}
\label{blocksystem}
\left[ \begin{array}{rr}
	A & E \\ 
	B & -I
\end{array}\right] 
\left[\begin{array}{r}
u \\ v
\end{array}\right] =
\left[\begin{array}{r}
f \\ 0
\end{array}\right].
\end{align}	
As is evident, any $[u,v]^T$ that solves this system also solves $Lu = f$, and satisfies $v = Bu$; this system represents splitting the range of $A$ and $B$ in obtaining a solution of $L$.   Taking one step of block Gaussian elimination of $A$ in this matrix to result in the upper triangular system
\begin{align}
\label{operator_decomp}
\left[ \begin{array}{rr}
	I & A^{-1} E \\ 
	0 & -I - BA^{-1}E
\end{array}\right] 
\left[\begin{array}{r}
u \\ v
\end{array}\right] =
\left[\begin{array}{r}
A^{-1}f \\ - BA^{-1}f
\end{array}\right],
\end{align}	
we then obtain $v$ as the solution to the following system, 
\begin{align}
	\label{schur_complement_system}
	(I + B A^{-1}E ) v = BA^{-1}f,
\end{align}
which represents the Schur complement system of $A$ in Eq.~(\ref{blocksystem}).  A back-substitution of $v$ into Eq.~(\ref{operator_decomp}),
\begin{align}
	u = A^{-1}(f - Ev),
\end{align}
results in $u$ the solution of $Lu = f$.  Because the SMPM is a discontinuous element discretization, $A$ is block diagonal and invertible, and so all divisions of $A$ are easily parallelized; the expensive part of the above is obtaining the solution of the Schur complement system.  In the rest of this paper, the focus is on efficiently obtaining this solution.  

\subsection{Inconsistency of the Poisson-Neumann system}
\label{inconsistency}
	Prior to obtaining the solution to Schur complement system, there remains the important point of dealing with the rank-deficiency of the Poisson-Neumann operator $L$.  The Poisson-Neumann equation is ill-posed in the continuous sense, and so the SMPM operator $L$ is rank-deficient and has non-trivial left and right null spaces of dimension one.  In symmetric discretizations, the kernel vector is the constant vector, but since $L$ is unsymmetric its left and right null spaces are different and only the right null space is constant vector.  To ensure consistency and solvability the right-hand-side vector $f$ is projected out of the left null space of $L$ \cite{Pozrikidis2001} and instead of $Lu = f$, the regularized system solved is
	\begin{align}
		\label{regularize_poisson}
		Lu =  \tilde{f}
	\end{align}
	where $\tilde{f} = f - u_L u_L^Tf$ is $f$ projected onto the range space of $L$ and $u_L \in \R^r$ is the unique vector with unit norm that satisfies $\norm{ u_L^T L }_2 = 0$.  The solution $u$ then is only known up to an indeterminant additive constant vector.  The rank deficiency of $L$ is inherited by the Schur complement system, and thus another regularization is required to project the Schur right hand side $b_S = BA^{-1}\tilde{f}$ out of the left null space of the Schur complement system.  Thus the Schur complement system,
	\begin{align}
		Sx_S = b_S,
	\end{align}
	 is modified to read
	\begin{align}
	\label{regularize_schur}
		Sx_S = b_S - u_S u_S^T b_S.
	\end{align}
	To summarize, the method for obtaining the solution $u$ to $Lu = f$ is shown in Algorithm \ref{alg:schur}.  The statement GMRES$(S,b_S)$ in Step \ref{schursolve} is meant to represent the solution of the linear system $Sx_S = b_S$ with the Generalized Minimum Residual Method (GMRES).
	
	\renewcommand{\algorithmicrequire}{\textbf{Input:}}
	\renewcommand{\algorithmicensure}{\textbf{Output:}}
	\begin{algorithm}[H]
		\begin{algorithmic}[1]	
			\REQUIRE $f, u_L, u_S$
			\ENSURE $u$
			\STATE $f \longleftarrow f - u_L u_L^T f$
			\STATE $b_S := BA^{-1}f$
			\STATE $b_S \longleftarrow b_S - u_S(u_S^T b_S)$  \label{reg2}
			\STATE $x_S :=  $ GMRES($S,b_S$) \label{schursolve}
			\STATE $u \longleftarrow A^{-1}(f - E x_S)$ \label{recover}
		\end{algorithmic}
		\caption{Schur complement method with null space projections.}
		\label{alg:schur}
	\end{algorithm}	

\section{Deflated Schur complement method}

	\begin{figure}[h]
	\begin{center}
		\includegraphics[width=0.35\textwidth]{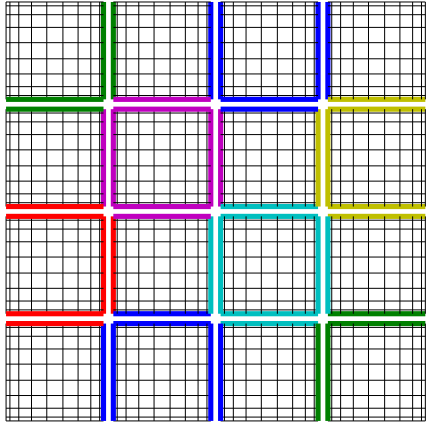}
	\end{center}
		\caption{ A sample domain with $m_x = 4$ and $m_y = 4$ elements in $x$ and $y$ respectively.  The elements have been separated to emphasize the discontinuous nature of the discretization.  The interfaces that together form the Schur grid $\Gamma$ are highlighted in color, with the colors signifying the blocks of the block-Jacobi preconditioner.  Interior blocks of the block-Jacobi preconditioner consist of eight interfaces; exterior blocks consist of six or four interfaces. }
		\label{fig_schur_preconditioner_grid}
	\end{figure}	

The iterative solution of the Schur complement system with GMRES requires an efficient preconditioner.  Many preconditioning techniques for the Schur complement system have been proposed \cite{Yamazaki2010a,Cros2002,Carvalho2001,Yamazaki2010} with most relying on two-level preconditioners: a local preconditioner that can be applied in parallel and a coarse global preconditioner to speed across-grid communication of components of the residual.
An example is the two-level additive-Schwarz preconditioner in which overlapping block-diagonal components are solved in parallel, augmented with a coarse grid correction to communicate information across the grid \cite{Fischer1997,Escobar-Vargas2014}).   In this work, a non-overlapping block-diagonal/block-Jacobi preconditioner is used, augmented with deflation, to achieve Krylov subspace convergence rates independent of the polynomial order, $p$, and weakly dependent on the element size, $h$.

First, note that for sparse matrices with non-zeros clustered around the diagonal, computing the inverse of  blocks along the diagonal separately can be a useful preconditioning technique.
Block-Jacobi preconditioners have been shown to be effective for the Schur complement of elliptic operators \cite{Couzy1995}, especially when combined with coarsened-grid preconditioners \cite{Manna2004, Pavarino2000, Pasquetti2006}.   Here, for preconditioning the Schur matrix, $S$, a block-Jacobi preconditioner is assembled in which a single block represents the coupling between the four interfaces bounding one element, and their corresponding interfaces in neighboring elements (for a total of $8n$ grid points in one block).  The elements in the $m_y \times m_x$ grid are divided in a checkerboard pattern as can be seen in Fig.~\ref{fig_schur_preconditioner_grid}, in which adjacent nodes are grouped together in blocks by color.  The colors correspond to the blocks of the block diagonal preconditioner for the Schur complement matrix.  Denoting the block-Jacobi preconditioner matrix as $M$, the preconditioned Schur complement system that is solved with GMRES is
\begin{align}
   SM^{-1}x'_S = b_S
\end{align}
and the solution is obtained by a final division by $M$
\begin{align}
   x_S = M^{-1}x'_S.
\end{align}
Since $M$ is explicitly block-diagonal (i.e. any non-zeros of $S$ coupling the blocks of $M$ are ignored in the factorization of $M$), divisions by $M$ can be computed efficiently in parallel.   The algorithmic summary of the preconditioned Schur complement method is given in Algorithm ~\ref{alg:pre}.

	\renewcommand{\algorithmicrequire}{\textbf{Input:}}
	\renewcommand{\algorithmicensure}{\textbf{Output:}}
	\begin{algorithm}[H]
		\begin{algorithmic}[1]	
			\REQUIRE $b, u_L, u_S, u_C$
			\ENSURE $x$
			\STATE $b \longleftarrow b - u_L u_L^T b$
			\STATE $b_S := BA^{-1}b$
			\STATE $b_S \longleftarrow b_S - u_S(u_S^T b_S)$ 
			\STATE $x :=  $ GMRES($SM^{-1},b_S$)
			\STATE $x \longleftarrow M^{-1} x$
			\STATE $x \longleftarrow A^{-1}(b - E x)$
		\end{algorithmic}
		\caption{Preconditioned Schur complement method}
		\label{alg:pre}
	\end{algorithm}
	
% The point: Brief summary of what deflation is.
Working in tandem with other preconditioners, deflation methods aim to accelerate the convergence of Krylov methods by eliminating (or ``deflating'') components of the residual within a chosen subspace.  The subspace is usually chosen to be a span of approximate eigenvectors of the operator corresponding to slowly converging eigenvalues.  Thus, the problematic eigenvalues are solved directly using a coarsened version of the operator, and the remaining components of the residual are eliminated by a Krylov solver.  Here, to augment the block-Jacobi preconditioner described in the previous section, a deflation method is used as the coarse-grid correction method, following the procedure in Ref.~\cite{Erlangga2008}.
  The deflation vectors are chosen to be a set of $d$ column vectors $Z \in \R^{k\times d}$ where $d \ll k$, and $k = \textrm{dim}(S)$.
   These deflation vectors are chosen to be discrete indicator vectors, equal to 1 on each pair of interfaces between two elements and zero everywhere else.
   Denoting as $\Gamma_j$ an interface between two elements, the $i$-th entry in the $j$-th deflation vector is given by
	\begin{align}
		(z_j)_i = \left\{
		     \begin{array}{ll}
		       1 & :  \textrm{ if $x_i \in \Gamma_j$ } \\
		       0 & : \textrm{ if $x_i \notin \Gamma_j$ }
		     \end{array}
		     \right\},
	\end{align}
	thus each vector is active on one pair of interfaces in the Schur grid.  The matrix of these vectors  $Z = [ z_1, z_2, \cdots, z_d ]$  defines a coarse version of the Schur problem, $C = Z^TSZ \in \R^{d\times d}$, and two projections
		\begin{align}
			P &= I - S ZC^{-1}Z^T \label{projection1}	 \\
			Q &= I - ZC^{-1}Z^TS
			\label{projection2}
		\end{align}
		each of size $\R^{k\times k}$.
      As a matrix $Z^T \in \R^{ d \times k}$ is a contraction operator that maps grid functions on the Schur grid to the coarse grid, and its transpose is a prolongation operator.
 The intuition behind the projections $P$ and $Q$ is that they project out of the subspace on which $ZC^{-1}Z^T$ is a good approximation of the left (in the case of $Q$) or right (in the case of $P$) inverse of $S$.  Thus the projections map onto the complement of the subspace on which the coarse matrix $C$ approximates the Schur matrix $S$ well.	  Finally, note that all applications of $C^{-1}$ require their own regularization since $C$ inherits rank deficiency from $S$;  denote as $u_C$ the left null space of $C$ in the following.  Deflation proceeds by noting that the solution of the Schur complement problem $Sx_S = b_S$ can be decomposed into
		\begin{align}
			x_S = (I - Q)x_S + Qx_S.
		\end{align}
		Then, the first term is just $ZC\backslash( Z^T - u_C u_C^TZ^T)(b_S - u_Su_S^Tb_S)$, which can be computed directly since $C$ is small. The second term is obtained by way of GMRES on the deflated and right-preconditioned system $PSM^{-1}x_S = P(b_S - u_Su_S^Tb_S)$ and then post-multiplying by $Q$, finally assembling the solution as
		\begin{align}
         \label{deflated_schur}
			x_S &= ZC\backslash( Z^T - u_C u_C^TZ^T)(b_S - u_Su_S^Tb_S)  \nonumber \\ &+ QM^{-1}\textrm{GMRES}(PSM^{-1},P(b_S - u_Su_S^Tb_S)).
		\end{align}
		Because $P$ projects out of the coarse space, the GMRES solution of $PSM^{-1}x_S = P( b_S - u_Su_S^Tb_S)$ minimizes only the component of the residual that cannot be well-approximated by the coarse solution.      This formulation of deflation-augmented right-preconditioning is an extension of the work in Ref.~\cite{Erlangga2008} to a rank-deficient matrix.  For completeness,  Algorithm \ref{alg:seq} depicts the algorithmic summary of the deflation method in which the notation GMRES$(A,b)$ is intended to represent the solution of a linear system $Ax = b$ with GMRES.
	\renewcommand{\algorithmicrequire}{\textbf{Input:}}
	\renewcommand{\algorithmicensure}{\textbf{Output:}}
	\begin{algorithm}[H]
		\begin{algorithmic}[1]	
			\REQUIRE $b, u_L, u_S, u_C$
			\ENSURE $x$
			\STATE $b \longleftarrow b - u_L u_L^T b$
			\STATE $b_S := BA^{-1}b$
			\STATE $b_S \longleftarrow b_S - u_S(u_S^T b_S)$ 
			\STATE $x_1 :=  $ GMRES($PSM^{-1},Pb_S$)
			\STATE $x_1 \longleftarrow QM^{-1} x_1$
			\STATE $x_2 := Z^T b_S - u_Cu_C^TZ^Tb_S$ 
			\STATE $x_2 \longleftarrow ZC\backslash x_2$
			\STATE $x := x_1 + x_2$
			%\STATE $x_C \longleftarrow x_C - v_C(v_C^Tx_C)$
			\STATE $x \longleftarrow A^{-1}(b - E x)$
		\end{algorithmic}
		\caption{Deflated and preconditioned Schur complement method}
		\label{alg:seq}
	\end{algorithm}

\section{Performance}

\begin{figure}[h]
	\begin{center}
		\includegraphics[width=0.95\textwidth]{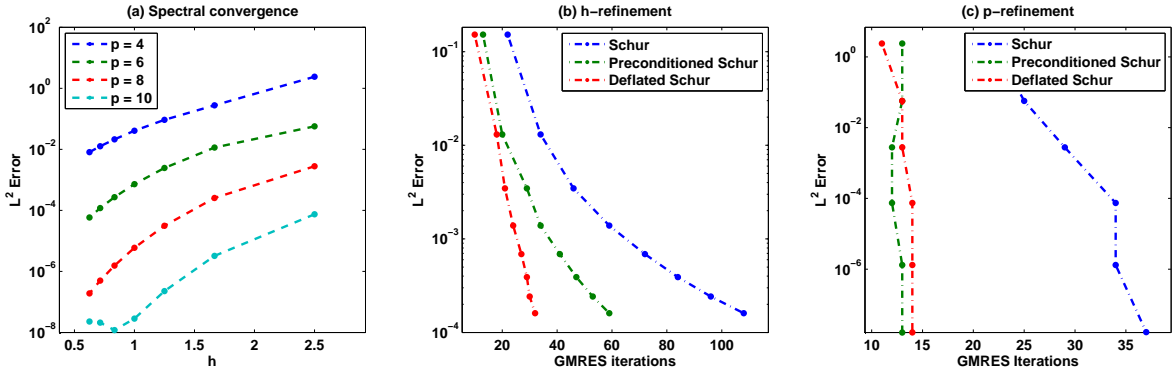}
	\end{center}
	\caption{\emph{Left}: The error decay as the initial grid with $p = m_x = m_y = 4$ is refined in both $h$ and $p$. Notice the decay of error is polynomial in $h$ but exponential in $p$.  \emph{Middle}:  $h$ refinement of the initial grid in which $p = m_x = m_y = 4$, with the analytic $L^2$ error plotted against the number of GMRES iterations required to achieve a relative tolerance of $10^{-10}$.  $h$ is gradually reduced as $m_x$ and $m_y$ are increased from $4$ to $32$. \emph{Right}:  $p$ refinement of the initial grid in which $p = m_x = m_y = 4$, with the analytic $L^2$ error plotted against the number of GMRES iterations required to achieve a relative tolerance of $10^{-10}$.  $p$ is gradually increased from $4$ to $24$. }
	\label{fig_results}
\end{figure}

To study the performance of the deflated Schur complement method, the domain $\Omega = [ 0, L_x ] \times [ 0, L_y ]$ is discretized with $n = 5$ (4th order polynomials), $m_x = 4$ and $m_y =4$.  On this grid the following Poisson problem was solved, 

	\begin{align}
		\nabla^2 u &= \cos( \lambda \pi x / L_x ) \cos( \lambda \pi y / L_y ) \nonumber \\
		n\cdot \nabla u &= 0,
	\end{align}
	with $\lambda \in \N$, which has the analytic solution
	\begin{align}
		u_a(x,y) = -\frac{L_x, L_y}{2 \lambda^2 \pi^2 }  \cos( \lambda \pi x / L_x ) \cos( \lambda \pi y / L_y ).
	\end{align}
	Besides evaluating the performance of the deflated Schur complement method with respect to GMRES convergence, comparing against an analytic solution makes it clear whether the discretization is exhibiting spectral convergence with respect to a known solution.  In the following subsections, we evaluate the GMRES convergence properties of the deflated Schur complement method as the initial grid is refined in both $p$, the polynomial order, and $h$, the element edge size ($L_x/m_x$).  The error $e$ is computed as the $L^2$ norm difference of the numerical solution, $u$, against the analytic solution, $u_a$, 
	\begin{align}
		e = \norm{u_a - u}_2.
	\end{align}
	In all of the examples below, $\lambda = 7$. 

	\subsection{Spectral convergence}

	First, note that as shown in Fig.~\ref{fig_results}(a),the Schur complement approach to solving the SMPM discretization converges to the true solution as a polynomial of $h$ and as an exponential of $p$.  This is made clear by the relatively gradual decay of the error in $h$ as compared with that in $p$, and is a hallmark of high-order methods like the spectral multi domain penalty method \cite{Escobar-Vargas2014}.  For the smallest values of $h$ and largest values of $p$ the error does not decay monotonically to zero; this may be an artifact of the ill-conditioning of the spectral differentiation matrices, which only worsens as $p$ grows and $h$ decreases.

	\subsection{$h$-refinement}
		
	Having established that the Schur complement method converges to the analytic solution at a rate equal to a polynomial of $h$, we now examine how the convergence of GMRES is affected as $h$ is refined.  Starting again with a grid $p = m_x = m_y = 4 $ on a domain $\Omega = [0, L_x] \times [ 0,  L_y]$, the number of elements $m_x$ and $m_y$ is iteratively grown, yielding a refinement in $h = L_x / m_x$.  The results of this refinement study are shown in Fig.~\ref{fig_results}(b) for values of $m_x, m_y = \{ 4, 8, 12, 16, 20, 24, 28, 32\}$. The Schur complement is assembled exactly as before, with decomposition along all internal element boundaries, and so the dimension of the Schur complement problem grows as $\mathcal{O}(m_xm_y)$.
	
	The GMRES algorithm is employed to reduce the residual to a relative tolerance of $10^{-10}$, at which point the rest of the Schur complement algorithm is employed to reconstruct the full solution (c.f. Step \ref{recover} in Algorithm \ref{alg:schur}).  The number of GMRES iterations required to achieve this tolerance is depicted along the horizontal axis in Fig.~\ref{fig_results}(b).  In the vertical axis of the same figure is shown the analytic error $\norm{u_h - u_a}_2$ in the resulting solution.  In all cases, the number of GMRES iterations grows as $h$ is refined.  This is first because the conditioning of the element stiffness matrices degrades as the element size goes to zero \cite{Costa2000}, and second because the dimension of the Schur complement matrix grows as $m_x$ and $m_y$ grow.
	
	However, deflating the Schur complement in the GMRES solver tempers the growth of the number of iterations as $h$ is refined.  At the finest grid, in which $m_x = m_y = 32$, the unpreconditioned Schur complement method takes nearly 110 Krylov iterations to converge; the deflated Schur complement method takes just over 30.  Furthermore, asymptotically, the number of iterations required to obtain a solution grows much more mildly in the deflated Schur complement approach.  Thus, it is observed that deflation, while not eliminating the dependence on $h$, strongly mitigates the growth in GMRES iterations as $h$ is reduced.  Finally, note that all of the Schur complement approaches shown in Fig.~\ref{fig_results}(b) are very efficient. For example, the number of grid points $r$ in the systems solved in Fig.~\ref{fig_results}(b) is $r = n^2 m_x m_y$.  Since $n=5$ and $m_x$ and $m_y$ grow to be 32 each, the total number of grid points grows to over $25,000$, which is greater by several orders of magnitude than the number of GMRES iterations required to obtain a solution that is correct to ten decimal places.

	\subsection{$p$-refinement}
	
	The results of a refinement study in $p$ are shown in Fig.~\ref{fig_results}(c) for values of $p = \{ 4, 6, 8, 10, 12, 14 \}$, and $m_x = m_y = 4$.  As $p$ is increased, the size of the Schur complement matrix grows as $p^2$, and its conditioning properties worsen due to the $h^{-p}$ conditioning of the spectral differentiation matrices embedded within it \cite{Costa2000}.  Nevertheless, it is observed in Fig.~\ref{fig_results}(c) that the convergence properties of GMRES are essentially unaffected by $p$-refinement when the block-diagonal preconditioner is applied to the Schur complement matrix.  Even the unpreconditioned Schur complement method only shows mild growth in GMRES iteration count as $p$ grows.  While it may seem that this $p$-independent convergence may depend on smoothness of the right-hand-side, this result has been confirmed for random white-noise right-hand-side vectors as well \cite{Joshi2016}.  It appears that GMRES convergence of the preconditioned Schur complement method is robust to refinements in $p$, which is particularly useful given that the error in the solution decays exponentially with $p$.  Finally, note again that the all variants of the Schur complement method shown in Fig.~\ref{fig_results}(c) are exceedingly efficient.  Since the grid grows as $p^2$, for the largest value of $p$, the number of grid points is $r = 3600$; yet, GMRES converges to ten digits of accuracy in less than 40 iterations in all cases, and in under 15 iterations in the preconditioned/deflated cases.
		
\section{Conclusion}

% Summary

A preconditioned Schur complement technique for solving the spectral multi-domain penalty method discretization of the Poisson-Neumann system was developed.  The preconditioning method relies on a local block-Jacobi preconditioner and subspace deflation to solve a coarse component of the residual.  By using both a local (block-Jacobi) and global (deflation) preconditioner, convergence of GMRES only mildly dependent on the grid resolution, $h$, and independent of polynomial order, $p$, is possible.  Since the error in the SMPM decays exponentially with $p$, achieving GMRES convergence independent of $p$ is very useful in practice, as it allows for high-accuracy solutions by increasing $p$ at minimal additional cost.   

% Impact: what are we using this in?
These ideas have already been leveraged within a high-order incompressible Navier-Stokes solver; the Poisson-Neumann system arises there as part of the incompressibility constraint on the flow.  By leveraging the deflated Schur complement algorithm described herein, problems on grids with millions of unknowns are solved with $\mathcal{O}(10)$ GMRES iterations.  The scalable efficient implementation of this deflated Schur complement approach is currently underway, but these methods show promise for parallelization since they are inspired from domain decomposition methods which are designed, in a sense, for distributed-memory parallel computing.  

% Future Work
A natural extension of this work is to generalize to three-dimensional problems.  We have shown already that if a periodic third dimension can be assumed, a Schur factorization can be used to efficiently extend these preconditioned Schur methods to three-dimensional problems \cite{Joshi2016}.  However, for more general three-dimensional approaches (in which periodicity cannot be assumed) more work is required to demonstrate the efficacy of the Schur complement method espoused in this work.  That having been said, there are  few theoretical hurdles for such an extension.  In particular, some Schur complement methods are plagued with so-called \emph{cross-points}, or grid points belonging to multiple elements, whose treatment requires special care in constructing the Schur complement problem.  Because the SMPM is a discontinuous discretization technique, there are no cross-points; every grid point belongs to only one element.  Thus, the extension of the Schur complement method to three dimensions is far less cumbersome in the SMPM than in continuous element methods. 

% Bibliography Stuff.
\bibliographystyle{plain}
\bibliography{capacitance_bibliography}

\end{document}